

\baselineskip=14pt
\parskip=10pt
\def\halmos{\hbox{\vrule height0.15cm width0.01cm\vbox{\hrule height
  0.01cm width0.2cm \vskip0.15cm \hrule height 0.01cm width0.2cm}\vrule
  height0.15cm width 0.01cm}}
\font\eightrm=cmr8 

\magnification=\magstephalf

\def\1{{\overline{1}}}
\def\2{{\overline{2}}}
\parindent=0pt
\overfullrule=0in

\def\frac#1#2{{#1 \over #2}}
\centerline
{\bf  D.H. Lehmer's Tridiagonal determinant:}
\centerline
{ \bf An \'Etude in (Andrews-Inspired) Experimental Mathematics }
\bigskip
\centerline
{\it Shalosh B. EKHAD and Doron ZEILBERGER}

\bigskip
\qquad\qquad\qquad\qquad \qquad\qquad\qquad\qquad  {\it Dedicated to George Andrews on his 80th birthday \quad } 

{\bf Lehmer's Theorem and its Finite Form}

Define, with Lehmer ([L], p. 54), $M(n)=M(n)(X,q)$, to be the following
tridiagonal $n \times n$ matrix (we changed $a$ to $\sqrt{X}$ and $r$ to $q$)
$$
M(n)_{i,j} \, = \,
\cases{ 1 \quad if \quad i-j=0 \, ; \cr 
 \sqrt{X} \, q^{(i-1)/2} \quad if \quad i-j=-1 \, ; \cr 
 \sqrt{X} \, q^{(i-2)/2} \quad if \quad i-j=1 \, ; \cr 
0 \quad otherwise \quad .}
$$

{\bf Theorem 1} (Lehmer [L])  
$$
\lim_{n \rightarrow \infty} \, \det \, M(n)(X,q) \, = \,
\sum_{a=0}^{\infty} \frac{(-1)^a X^a q^{a(a-1)}}{(1-q)(1-q^2) \cdots (1-q^a)} \quad .
$$

{\eightrm (As noted by Lehmer, when $X=-q$ and $X=-1$ one gets the sum sides of the famous Rogers-Ramanujan identities.)}

Our new result is an explicit expression for the {\bf finite form}, that immediately
implies Lehmer's theorem,  by taking the limit $n \rightarrow \infty$, and gives it a new (and shorter!) proof.

{\bf Theorem 2} 
$$
\det \, M(n)(X,q) \, = \,
\sum_{a=0}^{\lfloor n/2 \rfloor} \frac{(-1)^a X^a q^{a(a-1)} \, (1-q^{n-a}) (1-q^{n-a-1}) \cdots (1-q^{n-2a+1}) }{(1-q)(1-q^2) \cdots (1-q^a)} \quad .
$$

{\bf Proof of Theorem 2}: As noted by Lehmer ([L], Eq. (3)), by expanding with respect to the last row, we have
$$
\det \, M(n)(X,q) \, = \, \det \, M(n-1)(X,q) \, - \,  X\, q^{n-2} \, \det \, M(n-2)(X,q) \quad .
\eqno(LehmerRecurrence)
$$
Using the $q$-Zeilberger algorithm\footnote{$^1$}  
{\eightrm
Typing qzeil((-1)**a*X**a*q**(a*(a-1))*qbin(n-a,a),S,a,n,[]); in {\tt qEKHAD} gives
$X\,q^n-S+S^2, X\,q^n$ that is the recurrence operator annihilating the sum ($S$ is the forward shift operator in $n$) followed by the `certificate' (i.e. the proof, see [PWZ])) .}
([PWZ],[Z], see also [PR] for a nice Mathematica version),
we see that the right side of Theorem 2 also satisfies the very same recurrence. Since it holds for the {\bf initial conditions} $n=1$ and $n=2$ (check!),
the theorem follows by induction. \halmos

\vfill\eject

{\bf Secrets form the Kitchen}

Our paper could have ended here. We have increased human knowledge by extending a result of a famous number theorist,
and proved it rigorously.
{\it But}, at least as interesting as the {\bf statement} of the theorem (and far more interesting than the {\bf proof})
is the {\it way} it was {\it discovered}, and the rest of this paper will consist in describing {\bf two} ways of doing it.
The first way is a direct adaptation
of  George Andrews' ``{\it reverse-engineering}''  approach  
beautifully illustrated in the
last chapter of his delightful booklet [A]
({\eightrm based on 10 amazing lectures, given at Arizona State University, May 1985 
that we were fortunate to attend}). 
In that masterpiece (section 10.2)
he described how he used the computer algebra system SCRATCHPAD to prove a deep 
conjecture by three notable mathematicians: George Lusztig, Ian Macdonald and C.T.C. Wall.
In Andrews' approach, it is assumed that the discoverer knows about Gaussian polynomials, and
knows how to spot them. In other words the `atoms' are Gaussian polynomials. In the second, more
basic, approach, the only pre-requisite is the notion of {\it polynomials}, and Gaussian polynomials
pop-up naturally in the act of discovery.

We will start completely {\it from scratch}, pretending that we did not read Lehmer's paper. In fact
we did not have to `pretend'. We had no clue that Lehmer's paper existed until way after we discovered
(and proved) Theorem 2, (and hence reproved Lehmer's Theorem 1). This is the time for a short ``commercial break",
since this paper (like so many other ones!) owes it existence to the OEIS.                  

{\eightrm [Start of commercial break.]}

{\bf Serendipity and the OEIS}

We learned about Lehmer's Theorem 1 via {\it serendipity}, thanks to that amazing tool that we are so lucky to have, the
{\it On-Line Encyclopedia of Integer Sequences} [S] (OEIS).

Recall that a {\it composition} of $n$ is an array of positive integers $(p_1, \dots, p_k)$ such that
$p_1+ \dots+ p_k=n$, and they are very easy to count (there are $2^{n-1}$ of them). A {\it partition} of $n$
is a composition with the additional property that it is {\it weakly decreasing}, i.e.
$$
p_{i}-p_{i+1} \geq 0 \quad ( 1 \leq i < k \quad) ,
$$
(and they are much harder to count) .

My current PhD student, Mingjia Yang [Y] is investigating {\it relaxed partitions}, that she calls $r$-partitions, that
are compositions of $n$ with the condition
$$
p_{i}-p_{i+1} \geq r \quad .
$$
When $r=1$ we get the familiar partitions into distinct parts, and when $r=2$ we get one of the actors in the Rogers Ramanujan identities.
But what about negative $r$? In particular what about
$(-1)$-partitions? After generating the first twenty terms
$$
 1, 2, 4, 7, 13, 23, 41, 72, 127, 222, 388, 677, 1179, 2052, 3569, 6203, 10778, 18722, 32513, 56455 \quad ,
$$
we copied-and-pasted  it to the OEIS, and sure enough, we were {\it scooped}! It is sequence {\bf A003116},
 whose (former) description was `reciprocal of an expansion of  a determinant', that pointed to sequence {\bf A039924},
 mentioning  Lehmer's Theorem 1 (in fact the special case $X=q$).
As a reference it cited `personal communication' by  Herman P. Robinson,  a friend and disciple of Lehmer.
The OEIS entry for  {\bf A039924} also referenced  Lehmer's ``lecture notes on number theory'' 
but we could not find it either on-line or off-line.

Since Lehmer's proof seemed to have been lost, we tried to prove it ourselves, and succeeded.
Our approach, inspired by Andrews' [A], was to first find
an explicit expression for the {\it finite form}, and then take the limit as $n$ goes to infinity (like Andrews did
for the L-M-W conjecture). Only {\it after} we had the proof, we searched {\bf MathSciNet} for

``Lehmer AND determinant AND tridiagonal"  \quad ,

and discovered [L]. To our relief, Lehmer's proof was longer than ours, and
did not go via the finite form, Theorem 2. As far as we know, Theorem 2 is new.
Once we discovered the reference [L] we notified Neil Sloane, and he added that reference to the relevant sequences
{\bf A003116} and {\bf A039924}. So the present paper is yet another paper that owes its existence to the OEIS!

{\eightrm [End of commercial break.]}

{\bf How the Statement of Theorem 2 would have been (easily!) discovered by George Andrews}

In Andrews' proof of the L-M-W conjecture, he used the {\it Gaussian polynomials}  (aka as {\it q-binomial coefficients})
as {\bf building blocks}. With his approach, Theorem 2 could have been found by him fairly quickly.
Let $Q_n(X,q) := \det M(n)(X,q)$.

Recall that the Gaussian polynomials $GP(m,n)(q)$ are defined by
$$
GP(m,n)(q)  := \frac{(1-q^{m+1})  (1-q^{m+2}) \cdots  (1-q^{m+n})}{(1-q) \cdots (1-q^n)} \quad
$$
(in spite of their look, they are {\it polynomials!}) .

The way George Andrews would have discovered Theorem 2 is as follows.

Initially, crank out the first, say, twenty terms of the sequence of polynomials $Q_n(X,q)$, either by evaluating the determinants, or, more
efficiently,  via $(LehmerRecurrence)$.

You don't need a computer to realize that the coefficient of  $X^0$, i.e. the constant term, is always $1$.

The coefficients of $X=X^1$ in $Q_n(X,q)$ for $n$ from $1$ to $8$ are
$$
[0,-1,-1-q,-1-q-{q}^{2},-1-q-{q}^{2}-{q}^{3},-1-q-{q}^{2}-{q}^{3}-{q}^{4},-1-q-{q}^{2}-{q}^{3}-{q}^{4}-{q}^{5}, -1-q-{q}^{2}-{q}^{3}-{q}^{4}-{q}^{5}-q^{6}] \quad .
$$
A quick glance by George Andrews would have made him conjecture that it is
$$
-GP(n-2,1)(q) \quad .
$$

Moving right along, here are the coefficients of $X^2$ for $1 \leq n \leq 10$:

$$
[0,0,0,{q}^{2},{q}^{2}+{q}^{3}+{q}^{4},{q}^{2}+{q}^{3}+2\,{q}^{4}+{q}^{5}+{q}^{6},{q}^{2}+{q}^{3}+2\,{q}^{4}+2\,{q}^{5}+2\,{q}^{6}+{q}^{7}+{q}^{8},
$$
$$
{q}^{2}+{q}^{3}+2\,{q}^{4}+2\,{q}^{5}+3\,{q}^{6}+2\,{q}^{7}+2\,{q}^{8}+{q}^{10}+{q}^{9},{q}^{2}+{q}^{3}+2\,{q}^{4}+2\,{q}^{5}
+3\,{q}^{6}+3\,{q}^{7}+3\,{q}^{8}
+{q}^{11}+2\,{q}^{10}+2\,{q}^{9}+{q}^{12},
$$
$$
{q}^{2}+{q}^{3}+2\,{q}^{4}+2\,{q}^{5}+3\,{q}^{6}+3\,{q}^{7}+4\,{q}^{8}+2\,{q}^{11}+3\,{q}^{10}+3\,{q}^{9}+2\,{q}^{12}+{q
}^{13}+{q}^{14}] \quad .
$$
Dividing by $q^2$ and checking against the Gaussian polynomials `data base', suggests that the coefficient of $X^2$ is always
$$
q^2 \, GP(n-4,2)(q) \quad .
$$

Similarly, the coefficient of $X^3$ would have emerged as
$$
-q^6 \, GP(n-6,3)(q) \quad .
$$

The coefficient of $X^4$ would have emerged as
$$
q^{12} \, GP(n-8,4)(q) \quad .
$$

The coefficient of $X^5$ would have emerged as
$$
-q^{20} \, GP(n-10,5)(q) \quad .
$$

And {\it bingo}, it requires no great leap of an Andrews' imagination to conjecture that
$$
Q_n(X,q)=
\sum_{a=0}^{\lfloor n/2 \rfloor} (-1)^a \, X^a \, q^{a(a-1)} \, GP(n\,-\, 2a \, , \, a)(q) \quad ,
$$
that is identical to the statement of Theorem 2.

{\bf How the Statement of Theorem 2 could have been discovered by someone who is NOT George Andrews?}

Suppose that you have never heard of the Gaussian polynomials. You still could have  conjectured the statement of Theorem 2.
Even if you have never heard of {\it Gaussian} polynomials,  you probably {\it did} hear of
{\it polynomials}. So assuming the {\it ansatz} that , for each $a$, the coefficient of $X^a$ is a certain polynomial
in $q^n$, try and fit it with a `generic' polynomial with {\it undetermined coefficients}.\footnote{$^2$}
{\eightrm You start out with a generic polynomial of degree $0$, and keep raising the degree until success (or failure).}

Setting $N=q^n$, your computer would have guessed the following polynomial expressions (in $N=q^n$) for the first five coefficients of $X$ in $Q_n(X)$.

$\bullet$ The coefficient of $X$ in $Q_n(X,q)$ is
$$
{\frac {N-q}{q \left( 1-q \right) }} \quad .
$$

$\bullet$ The coefficient of $X^2$ in $Q_n(X,q)$ is
$$
{\frac { \left( N-{q}^{2} \right)  \left( N-{q}^{3} \right) }{{q}^{3} \left( 1+q \right)  \left( 1- q \right) ^{2}}} \quad .
$$

$\bullet$ The coefficient of $X^3$ in $Q_n(X,q)$ is
$$
-{\frac { \left( N-{q}^{3} \right)  \left( N-{q}^{4} \right)  \left( N-{q}^{5} \right) }{{q}^{6} \left( 1+q \right)  \left( {q}^{2}+q+1 \right)  \left( q-1
 \right) ^{3}}} \quad .
$$

$\bullet$ The coefficient of $X^4$ in $Q_n(X,q)$ is
$$
{\frac { \left( N-{q}^{4} \right)  \left( N-{q}^{5} \right)  \left( N-{q}^{6} \right)  \left( N-{q}^{7} \right) }{{q}^{10} \left( {q}^{2}+1 \right)  \left( 
q-1 \right) ^{4} \left( 1+q \right) ^{2} \left( {q}^{2}+q+1 \right) }} \quad .
$$

$\bullet$ The coefficient of $X^5$ in $Q_n(X,q)$ is
$$
-{\frac { \left( N-{q}^{5} \right)  \left( N-{q}^{6} \right)  \left( N-{q}^{7} \right)  \left( N-{q}^{8} \right)  \left( N-{q}^{9} \right) }{{q}^{15}
 \left( q-1 \right) ^{5} \left( {q}^{4}+{q}^{3}+{q}^{2}+q+1 \right)  \left( 1+q \right) ^{2} \left( {q}^{2}+q+1 \right)  \left( {q}^{2}+1 \right) }} \quad .
$$

This immediately leads one to guess that the {\bf numerator} is always
$$
(-1)^a \, (N-q^a)(N-q^{a+1}) \cdots (N-q^{2a-1}) \quad .
$$
On the other hand, the sequence of denominators, let's them call them $d(a)$, for $1 \leq a \leq 5$, happens to be
$$
[-q \left( q-1 \right) ,{q}^{3} \left( 1+q \right)  \left( q-1 \right) ^{2},-{q}^{6} \left( 1+q \right)  \left( {q}^{2}+q+1 \right)  \left( q-1 \right) ^{3}
,{q}^{10} \left( {q}^{2}+1 \right)  \left( q-1 \right) ^{4} \left( 1+q \right) ^{2} \left( {q}^{2}+q+1 \right) 
$$
$$
,-{q}^{15} \left( q-1 \right) ^{5} \left( {q}
^{4}+{q}^{3}+{q}^{2}+q+1 \right)  \left( 1+q \right) ^{2} \left( {q}^{2}+q+1 \right)  \left( {q}^{2}+1 \right) ] \quad .
$$
This looks a  bit complicated, but let's form the sequence of {\it ratios} $d(a)/d(a-1)$ for $a=2,3,4,5$  and expand, getting
$$
[{q}^{2}-{q}^{4},{q}^{3}-{q}^{6},-{q}^{8}+{q}^{4},{q}^{5}-{q}^{10}] \quad,
$$
that is clearly $q^a (1-q^a)$. Hence the coefficient of $X^a$ in $Q_n(X,q)$ is guessed to be
$$
\frac{(-1)^a \, (N-q^a)(N-q^{a+1}) \cdots (N-q^{2a-1})}{q^{a(a+1)/2} \, (1-q) \cdots (1-q^a)} \quad .
$$
By putting $N=q^n$ we get the statement of Theorem 2.

So with  this second approach, we discovered the Gaussian polynomials {\it ab initio}, our only gamble was that the coefficients of $X$ in $Q_n(X,q)$
are always {\it polynomials} in $q^n$.

{\bf Concluding words}

Let us quote the last sentence of section 10.2 of [A], where Andrews described his pioneering ({\it experimental mathematics!}) approach 
illustrated by his discovery process of the proof of the L-M-W conjecture.

{\it ``From here the battle with the L-M-W  conjecture is 90 percent won. Standard techniques allow one to establish
the [finite form] of the conjecture, and a simple argument leads to the original conjecture.''}

Today the $90$ percent may be replaced by $99.999$ percent, since the final {\it verification} can be done
automatically by using the so-called $q$-Zeilberger algorithm ([PWZ][Z][PR]). In the much more difficult L-M-N case, this would have
saved George Andrews a few hours, and would have made it accessible to anyone else. 
In the present case, you can still use 
the $q$-Zeilberger algorithm, if you are feeling lazy, but it is not too hard to do it purely humanly.
Can {\bf you} do it?

{\bf References}

[A] George E. Andrews, {\it ``q-Series: Their Development and Applications in Analysis, Number Theory, Combinatorics, Physics, and
Computer Algebra''}, CBMS \#66, American Mathematical Society, 1986.

[L] Derrick H. Lehmer, {\it Combinatorial and cyclotomic properties of certain tridiagonal matrices}, Proceedings of the Fifth Southeastern
Conference on Combinatorics, Graph Theory, and Computing (Florida Atlantic University, Boca Raton, Fla. 1974)), 53-74.
Congressus Numerantium, No. X, Utilitas Math., Winnipeg., Man. 1974. {\bf MR0441852(56 \#243)} \quad . 
A scanned copy is available from \hfill\break
{\tt http://sites.math.rutgers.edu/\~{}zeilberg/akherim/LehmerDet1974.pdf } \quad  [Accessed Aug. 18, 2018].

[PR] Peter Paule and Axel Riese, {\it A Mathematica q-analogue of Zeilberger's algorithm based on an algebraically motivated approach to $q$-hypergeometric Telescoping}, 
in: ``Special Functions, q-Series and Related Topics,'' Fields Inst. Commun. {\bf 14}, 179-210, 1997 \quad .  \hfill\break Available from:
{\tt http://www.risc.jku.at/publications/download/risc\_118/diss.pdf } \quad  ~[Accessed Aug. 18, 2018].

[PWZ] Marko Petkov{\v s}ek, Herbert S. Wilf, and Doron Zeilberger, {\it ``A=B"}, A.K. Peters, 1996. \hfill\break Available from:
{\tt https://www.math.upenn.edu/\~{}wilf/Downld.html} ~[Accessed Aug. 18, 2018].\hfill

[S] Neil A. J. Sloane, {\it The On-line Encyclopedia of Integer Sequences} \quad , \quad {\tt  https://oeis.org} \quad .

[Y] Mingjia Yang, {\it Relaxed partitions}, in preparation.

[Z] Doron Zeilberger, {\tt qEKHAD} (a Maple program implementing the $q$-Zeilberger algorithm). Available from \hfill\break
{\tt http://sites.math.rutgers.edu/\~{}zeilberg/tokhniot/qEKHAD }  \quad [Accessed Aug. 18, 2018].

\bigskip
\hrule
\bigskip
Shalosh B. Ekhad, c/o D. Zeilberger, Department of Mathematics, Rutgers University (New Brunswick), Hill Center-Busch Campus, 110 Frelinghuysen
Rd., Piscataway, NJ 08854-8019, USA. \hfill\break
Email: {\tt ShaloshBEkhad at gmail dot com}   \quad .
\bigskip
Doron Zeilberger, Department of Mathematics, Rutgers University (New Brunswick), Hill Center-Busch Campus, 110 Frelinghuysen
Rd., Piscataway, NJ 08854-8019, USA. \hfill\break
Email: {\tt DoronZeil at gmail  dot com}   \quad .

\end